\documentclass[12pt]{article}
\usepackage{amsmath}
\usepackage{latexsym}
\usepackage{amsfonts}
\usepackage[normalem]{ulem}
\usepackage{array}
\usepackage{amssymb}
\usepackage{graphicx}
\usepackage{mathtools}
\usepackage{transparent}
\usepackage{amsmath}
\usepackage{xcolor}
\usepackage[skip=2ex]{caption}
    \usepackage[belowskip=2ex]{subcaption}
    \usepackage{lipsum}
    \usepackage{graphicx}
        \setkeys{Gin}{width=.55\textwidth}
    \usepackage{titlesec}
    \usepackage{xcolor}
    \usepackage{ulem}
    \usepackage{amsmath}
    \usepackage{xfrac}
    \usepackage{amssymb}
    \usepackage{amsfonts}
    \usepackage{physics}
    \usepackage{slashed}
    
\usepackage{epigraph} %Para la frase del principio
\usepackage[toc,page]{appendix}
\usepackage[backend=biber,
style=numeric,
sorting=none,
isbn=false,
doi=false,
url=false,
]{biblatex}

\usepackage{wrapfig}
\usepackage{wasysym}
\usepackage{enumitem}
\usepackage{adjustbox}
\usepackage{ragged2e}
\usepackage{nicefrac} %para las fracciones fancy
\usepackage{longtable}
\usepackage{changepage}
\usepackage{setspace}
\usepackage{hhline}
\usepackage{multicol}
\usepackage{tabto}
\usepackage{float}
\usepackage{multirow}
\usepackage{makecell}
\usepackage{fancyhdr}
\usepackage[toc,page]{appendix}
\usepackage[hidelinks]{hyperref}

\usepackage{flowchart}\usepackage[paperheight=11.69in,paperwidth=8.27in,left=0.98in,right=0.99in,top=1.08in,bottom=0.79in,headheight=1in]{geometry}
\usepackage[utf8]{inputenc}
\usepackage[T1]{fontenc}
\TabPositions{0.5in,1.0in,1.5in,2.0in,2.5in,3.0in,3.5in,4.0in,4.5in,5.0in,5.5in,6.0in,}
\usepackage{titling} %Esto es para poner el titulo con abstract
    \title{A heuristic study of prime number distribution}
    \author{Carlos Ros Pérez}
    \date{\today}
\usepackage{fancyhdr} %Esto controla el header
\fancypagestyle{plain}{}
\begin{document}

\pagestyle{plain} %Numeraria las páginas si no usaras fancyhdr
\setlength{\parskip}{4.0pt} %Mejora el interlineado
\setlength\parindent{24pt}  %%Esto controla el indent

\maketitle
\thispagestyle{empty}
\begin{abstract}
This work consists of a heuristic study on the distribution of prime numbers in short intervals. We have modelled the occurrence of prime numbers such intervals as a counting experiment, as a result, we have provided an experimental validation and an extension of the Montgomery and Soundararajan conjecture. This is a reduced version of my bachelor's thesis presented at the University of Valencia.
\end{abstract}

\section{Introduction}
\begin{flushright}
\rightskip=1.8cm\textit{“It is evident that the primes are randomly distributed but,\break unfortunately, we don’t know what ‘random’ means.”} \\
\vspace{.2em}
\rightskip=.8cm---R. C. Vaughan
\end{flushright}

\indent At the end of the eighteenth century, at the age of 15, Gauss empirically found that the density of primes in the neighborhood of an integer n can be estimated by $\frac{1}{\log n}$ \cite{1}, leading him to conjecture the celebrated prime number theorem\footnote{$f(x) \sim g(x)$ means$ f(x)$ is asymptotically equivalent to $g(x)$,i.e., $f(x) \sim g(x) 
\rightarrow \underset{x \to\infty}{\lim}\frac{f(x)}{g(x)}=1 $}:

\begin{equation}
   \pi(x)\sim \frac{x}{\log x}
    \label{eq:PNT}
\end{equation}

This conjecture  was demonstrated independently and at the same time by Hadamard \cite{2} and de la Valée Poussin \cite{3}.
%In the following, we are going to use the notation $li(x)= \int_{0}^{x}\frac{dt}{\log t}$\footnote{The integral $li(x)$ has to be interpreted as a Cauchy principal value.} and $Li(x)= \int_{2}^{x}\frac{dt}{\log t}$. 

\indent In the first half of the twentieth century, Cramér tried a revolutionary approach consists in considering the prime numbers as random variables, considering that, for a given a number N, there is a chance of $\frac{1}{\log N}$ of being prime, therefore a probability of $1-\frac{1}{\log N}$  of being composed, and the different numbers are treated as independent events. This kind of approach has many shortcomings, such as considering that the probability of one number and the following are primes is different from zero, or provide a non-vanishing probability that an even number is prime. Nevertheless, it has had many successes, for instance leading to the asymptotic limit: $\underset{p_{n} \leq x}{\max}(p_{n+1}-p_{n}) \sim \log ^{2} x $, known as the Crámer conjecture \cite{4}. That conjecture can be supplemented by the Hardy–Littlewood Prime k-tuples Conjecture. To understand what it is about, it is necessary to introduce the concept of k-tuple. We may define a k-tuple ($\wp _{k}$) simply as a set  of k distinct non-negative integers. We can define the following function:
\begin{equation}
   \mathcal{L}(\wp _{k})= \prod_{p} \frac{p^{k-1}(p-\nu_{\wp _{k}}(p))}{(p-1)^{k}}
    \label{eq:HLconstant}
\end{equation}

\indent Where the product is extended over all prime numbers and where $\nu_{\wp _{k}}(p)$ is the number of distinct residue classes modulo p represented by elements of $\wp _{k}$ \footnote{ In other words, if we consider the list formed by all the remainders left by dividing the different numbers in $\wp _{k}$ by p, the number residue classes modulo p occupied by $\wp _{k}$ would correspond to the number of different elements in the former list}. Let   $\wp _{k}=\left \{h_{1},h_{2},...,h_{k} \right \}$, we denote as $\pi(x;\wp _{k})$ the number of positive integers $n \leq x$ such that $\left \{n+ h_{1},n+ h_{2},...,n+ h_{k}\right \}$ are prime numbers. Now we can enunciate the Hardy–Littlewood k-tuples conjecture \cite{5}:
\begin{equation}
   \pi(x;\wp _{k}) \sim \mathcal{L}(\wp _{k}) \frac{x}{\log^{k} x}
    \label{eq:HLktuple}
\end{equation}

Now the definition of $\nu_{\wp _{k}}$ makes sense, because if the elements of $\wp _{k}$ occupy all the residue classes modulo p, resulting in $\mathcal{L}(\wp _{k})=0$, then at least one number of the list $\left \{n+ h_{1},n+ h_{2},...,n+ h_{k}\right \}$  must be divisible by p, forcing $\pi(x;\wp _{k})$ to be zero.

\indent In 1975 Gallagher, starting from the Hardy–Littlewood k-tuples Conjecture, demonstrated the following result \cite{6}:
\begin{equation}
   \# \left \{integers\; x \leq N \;/\; \pi(x+\lambda \log x)-\pi(x)=k\; \right \} \sim N \lambda^{k} \frac{e^{-\lambda}}{k!}
    \label{eq:poisson}
\end{equation}

Which is a remarkable result because it states that the distribution of the number of integers (x) such that there are k prime numbers between x and  x+$\lambda \log x$ asymptotically tends to the well-known Poisson distribution. This result can be generalized to assess the distribution of primes in some intervals. Mongomery and Soundararajan \cite{7} demonstrated\footnote{In the article they work with the second Chebyshev function $\psi(x)$ instead of with $\pi(x)$, we can relate its moments in the following way:  $M_{k}(\pi,N)=\frac{M_{k}(\psi,N)}{(\log N)^{k}}$} that the Cramér model correctly predicts the distribution of primes in intervals of length $h(N)$ where $h \asymp \log N$ \footnote{$f(x) \asymp g(x) \rightarrow f(x)=\mathcal{O}(g(x))$ and $g(x)=\mathcal{O}(f(x)) $} with  mean and variance $\sim \frac{h}{\log N}$, this range is called "microscopic" scale \cite{8}. Nevertheless, it fails out of that range, which led them to conjecture that if $(\log N)^{1+\delta} \leq h \leq (N^{1-\delta}$)\footnote{In particular, we can set $\delta\rightarrow \infty$, then the condition becomes $\frac{h}{\log N} \sim \infty$ and $\frac{h}{N} \sim 0$}, what is called "mesoscopic" scale, then the distribution is normal with mean $\sim \frac{h}{\log N}$ and variance $\sim \frac{h}{\log^{2} N} \log \frac{N}{h}$.

\indent In this work we are going to study the distribution of the prime numbers on short scales with the aid of the methods developed by Sanchis \cite{9}.

\section{Computational model}

\hspace{\parindent} The main purpose of this work is to conduct a counting experiment on prime numbers in short intervals to the study their distribution. In particular, we are going to compute the mean $\langle p \rangle$ and variance $\sigma_{p}$ of the number of primes in certain intervals, supplying an experimental test of the predictions of Mongomery and Soundararajan (hereinafter cited as MS). Furthermore, we are going to look for an empirical formula for $\sigma_{p}$ in intervals with finite N.

\begin{figure}[!t]
\centering
\captionsetup{justification=centering}
\includegraphics[scale=1]{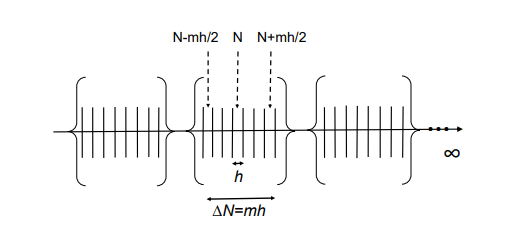}
\captionsetup{width=0.95\linewidth}
            \captionsetup{justification=centering}
\caption{Interval around an integer number N (taken from \cite{9}) }
\label{fig:2}
\end{figure}

\indent We are going to work with number intervals around an integer N. These intervals are divided into m subsets of length h, see figure \ref{fig:2} for a diagrammatic view. Let p be a random variable symbolising the number of prime numbers in a subset. Given intervals with fixed h and m, we are going to study the normalized variance $w$ as we vary N.
\begin{equation}
   w=\frac{\sigma^{2}_{p}}{\langle p \rangle}=\frac{\langle p^{2} \rangle-\langle p \rangle^{2}}{\langle p \rangle}
   \label{eq:norvar}
\end{equation}

\indent With $\langle p \rangle=\frac{\sum_{k=1}^{m}p_{k}}{m}$, where $p_{k}$ is the number of primes in the k-th subset. In order to parametrize $w$ we must recall the characteristic scales \cite{8}:
\begin{enumerate}
  \item Microscopic scale: When $h \asymp \log N$. On this scale the mean and the variance tend asympotically to  $\frac{h}{\log N}$
  \item Mesoscopic scale: When $\frac{h}{\log N} \rightarrow \infty$ and $\frac{h}{N} \rightarrow 0$. On this scale $\langle p \rangle \sim \frac{h}{\log N}$ and $\sigma_{p}^{2}\sim \frac{h}{\log^{2} N} \log \frac{N}{h}$.
  \item Macroscopic scale: When $h(N)>>>N$. On this scale there is no known distribution law.
\end{enumerate}

% The definition of microscopic scale is an asymptotic limit, however, for finite h, we may consider that the microscopic scale is reached when

In this work, we are going to employ fixed values of h such that $h<<<N$ \footnote{Quantitatively, the biggest h that we are employing is $h=5\cdot 10^{4}$ and the smallest N is $N=10^{9}$}, therefore we never reach the macroscopic scale. The definitions of the scales are asymptotic limits, however, for finite and fixed h, we may consider that an interval is mesoscopic if $\frac {h} {\log N} >> 1 $ is verified, otherwise, we would consider it microscopic. We shall take into account several facts. For fixed h, $N \rightarrow \infty$ leads to $\frac {h} {\log N} \rightarrow 0$, hence, we expect a Poisson distribution in that limit, thus $\underset{N \rightarrow \infty}{\lim}w=1$. In addition, for large values of h, we may consider a value for the mean and the variance similar to the predicted by MS in the asymptotic limit ($\frac{h}{\log N}$ and $\frac{h}{\log^{2} N} \log \frac{N}{h}$ respectively), therefore $w \sim 1-\frac{\log h}{\log N}$. We can try to model the normalized variance as:
\begin{equation}
   w=1-\frac{b(h,m)}{\log N}
   \label{eq:norvarg}
\end{equation}

\indent It is noteworthy to mention that we are not going to propose a model for the mean of the distribution because, as we will see to see in our results, the asymptotic value $\langle p \rangle \sim \frac{h}{\log N}$ predicts it with an astonishing accuracy, even for N as small as $10^{9}$. Nevertheless, if we had chosen $w=1-\frac{\log h}{\log N}$, it would have led to inaccurate results. The reason is that the MS correction is not providing a formula to compute the variance and the mean for large h, but an asymptotic limit for those quantities. Thus, we are measuring the mean and variance of an unknown distribution and expecting, based on the MS correction, that they are similar to $\frac{h}{\log N}$ and $ \frac{h}{\log^{2} N} \log \frac{N}{h}$ respectively, in our range of $h$ and $N$. This means that for large N and h we may assume $b(h,m)=\log h$. In this context, we consider h large if $\frac {h} {\log N} >>> 1 $.
%%Aqui no estamos cogiendo <p> directamente como el <p> de la escala mesoscopica o microscopica aunque en los fits nos salga que el mismo para ambas. Eso es debido a que <p> se da para límites asintoticos, por lo tanto consideramos una desviación de esta.
%IMPORTANTE se afirma que si  we are neither on the microscopic scale nor on the mesoscopic one. Esto es porque cuando escribi el texto no acababa de entender el planteamiento, tal y como estamos considerando las escalas si estamos en la escala microscópica. Lo podre justificar diceindo que las escalas se refieren a una h como función de N, como no estamos considerando h(N) si no una h fija, no podemos decdir que estamos en ninguna escala. Problema, la escala mesoscopica viene de considerar (log(N))^(1+d)<h<n^(1-d), por lo tanto si consideramos h=1 necesariamente no se cumple dicha condición

%Recuerda que el hecho de incluir la constante C es pora obtener una formula lo más general posible, ya que sabemos que al tender al mismo valor de la media en ambas escalas esa C ha de ser 1

\indent On the other hand, a subinterval of length $h=1$ may contain either one prime or none at all, implying $p_{k}\in \{0,1\}$, and therefore $\langle p \rangle=\frac{\sum_{k=1}^{m}p_{k}}{m}=\frac{\sum_{k=1}^{m}p_{k}^{2}}{m}= \langle p^{2} \rangle $. Consequently, for $h=1$, $w=1-\langle p \rangle$. If we choose $h=1$ we may consider that we are at the microscopic scale, therefore, a suitable choice of $\langle p \rangle$ would be a deviation from $\frac{1}{\log N}$, we may parametrize it as $\langle p \rangle=\frac{C}{\log N}+D$,with C an arbitrary constant that has been include to extend the generality of the formula, but we expect, and the result shall confirm, that $C=1$. On the other hand, a Poisson distribution is expected in the limit $N \rightarrow  \infty$, which implies $w=1$, hence $D=0$.

\indent To sum it all up, we expect $b(1,m)=C$ and $b(h,m)=\log h$ for large h. If we combine both limits, the simplest estimate for $b(h,m)$ could be $b(h,m)=C+\alpha (h) \log h$, with $\underset{h \rightarrow  \infty}{\lim} \alpha (h)=1 $. The expression of $\alpha (h)$ will be estimated in base of the empirical results. Therefore, the normalized variance would read:

\begin{equation}
   w=1-\frac{C+\alpha (h) \log{h}}{\log{N}}
   \label{eq:13}
\end{equation}

\indent As in \cite{9} we are going to make the simulations with intervals with three different values of m: $m=10^{3}$, $m=10^{4}$ and $m=10^{5}$, using different values of h for every m, and varying N for every pair $(h,m)$. 
\subsection{Error analysis}

\hspace{20.0pt} As every empirical procedure, our calculations are not exempt from errors. In this section we aim to estimate the error associated with the number of primes in each interval. In order to assess the value of the errors, we are going to consider the probability distribution given by the asymptotic limits of the MS conjecture, i.e., a normal distribution characterized by $\mu=\frac{h}{\log N}$ and $\sigma^{2}=\frac{h}{(\log N)^{2}}\log \frac{N}{h}$. In this work we are computing the statistical estimators from the value of the random variable $p$ in the different subsets and associating the result to an interval of fixed N. Nevertheless, we are not considering the N associated to each subset but the N associated to the interval. Therefore we are inferring a systematic error in the results. Considering that the length of the interval is $\Delta N=m \cdot h$. It is possible to assess the value of this error in the normal approximation as the difference between the maximum and the minimum value of the mean of primes in the interval.
\begin{equation}
   \Delta \mu=\frac{h}{\log (N-\frac{\Delta N}{2})}-\frac{h}{\log (N+\frac{\Delta N}{2})}=\frac{2h \tanh^{-1}(\frac{\Delta N}{2 N})}{\log (N-\frac{\Delta N}{2}) \log (N+\frac{\Delta N}{2})}
   \label{eq:siserr}
\end{equation}

\indent Where we have used $\tanh^{-1}(x)=\frac{1}{2}\log(\frac{1+x}{1-x})$. By Taylor expanding equation (\ref{eq:siserr}) around $N \rightarrow \infty$, keeping the first-order terms and substituting $\Delta N=m\cdot h$ we obtain 
\begin{equation}
   \Delta \mu=\frac{mh^{2}}{N (\log N)^{2}}
   \label{eq:siserr2}
\end{equation}
\indent We are led to the following relative systematic error:

\begin{equation}
   \epsilon_{sys}=\frac{\Delta \mu}{\mu}=\frac{mh}{N \log N}
   \label{eq:siserr3}
\end{equation}

\indent We can obtain the systematic error of the variance using the same procedure, by considering the value of the variance in the asymptotic limit we obtain:

\begin{equation}
    \Delta \sigma_{p}^{2}=\frac{h\log(\nicefrac{\left (N-\frac{\Delta N}{2} \right )}{h})}{\log (N-\frac{\Delta N}{2})^{2}}-\frac{h\log(\nicefrac{\left (N+\frac{\Delta N}{2} \right )}{h})}{\log (N+\frac{\Delta N}{2})^{2}}\approx \frac{h \Delta N \left (2\log(\nicefrac{N}{h})-\log(N) \right )}{\log(N)^{3}N}
   \label{eq:siserr3}
\end{equation}

\indent We can easily compute the relative statistical error of $\langle p \rangle$ given that we are doing a counting experiment related to a normal random variable with standard deviation\break
$\sigma=\sqrt{\frac{h}{(\log N)^{2}}\log \frac{N}{h}}$.
%%No tenemos en cuenta el valor de m en el calculo porque estamos calculando la varianza de la función <p>=\sum p_{k}/N, por propagacion de errores su varianza es N \sigma_{p} /N =\sigma_{p}
\begin{equation}
   \epsilon_{stat}=\frac{\sigma(p) }{\mu}=\sqrt{\frac{1}{h}\log\frac{N}{h}}
   \label{eq:staterr}
\end{equation}

\indent The range in which we are interested is $\frac{h}{\log N} >>> 1$ and $\frac{h}{N} <<< 1$, therefore it is a suitable scope to keep errors small. Since we are restricted to the computing power of an average computer, we cannot approach these limits, so errors are not negligible. Quantitatively; for $m=10^{5}$, $h=10^{3}$ and $N=10^{10}$; the relative systematic and statistic errors are $\epsilon_{sys} =0.04\%$ and $\epsilon_{stat} =12\%$. 

\indent The parameter m is related to the number of measures of primes carried out, hence, the bigger m, the more accurate the data will be. While keeping a low systematical error.

\begin{table}[b]
\label{tab:17}
\resizebox{\textwidth}{!}{%
\begin{tabular}{|l|c|c|c|c|c|c|}
\hline
N                                                                         & $10^{9}$ & $10^{10}$ & $10^{11}$ & $10^{12}$ & $10^{13}$ & $10^{14}$ \\ \hline
Empirical mean                                                            & $120,655 \pm 1.455 \pm 8.665$  & $108,568 \pm 0.118 \pm 8.466$   & $98,724 \pm 0.009 \pm 8.259$    & $90,483 \pm 8.053$    & $83,517 \pm 7.854$    & $77,536 \pm 7.663$    \\ \hline
\begin{tabular}[c]{@{}l@{}}Theoretical asymptotic\\ mean\end{tabular}     & $120,637$  & $108,574$   & $98,703$    & $90,478$    & $83,518$    & $77,553$    \\ \hline
Empirical variance                                                        & $67,380\pm 0.356$   & $65,103\pm 0.037$    & $62,712\pm 0.004$   & $60,048$    & $57,531$    & $55,300$    \\ \hline
\begin{tabular}[c]{@{}l@{}}Theoretical asymptotic\\ variance\end{tabular} & $75,091$   & $71,681$    & $68,213$    & $64,858$    & $61,688$    & $58,730$    \\ \hline
\end{tabular}}

\captionsetup{width=0.95\linewidth}
            \captionsetup{justification=centering}
            
\caption{ Comparison between empirical and theoretical parameters for an interval with $h=2500$ and $m=10^{5}$ }
\end{table}

 \begin{figure}[!b]\centering
        \setlength\belowcaptionskip{-1.5ex}
          \begin{subfigure}[t]{0.45\textwidth}
           \includegraphics[width=\textwidth]{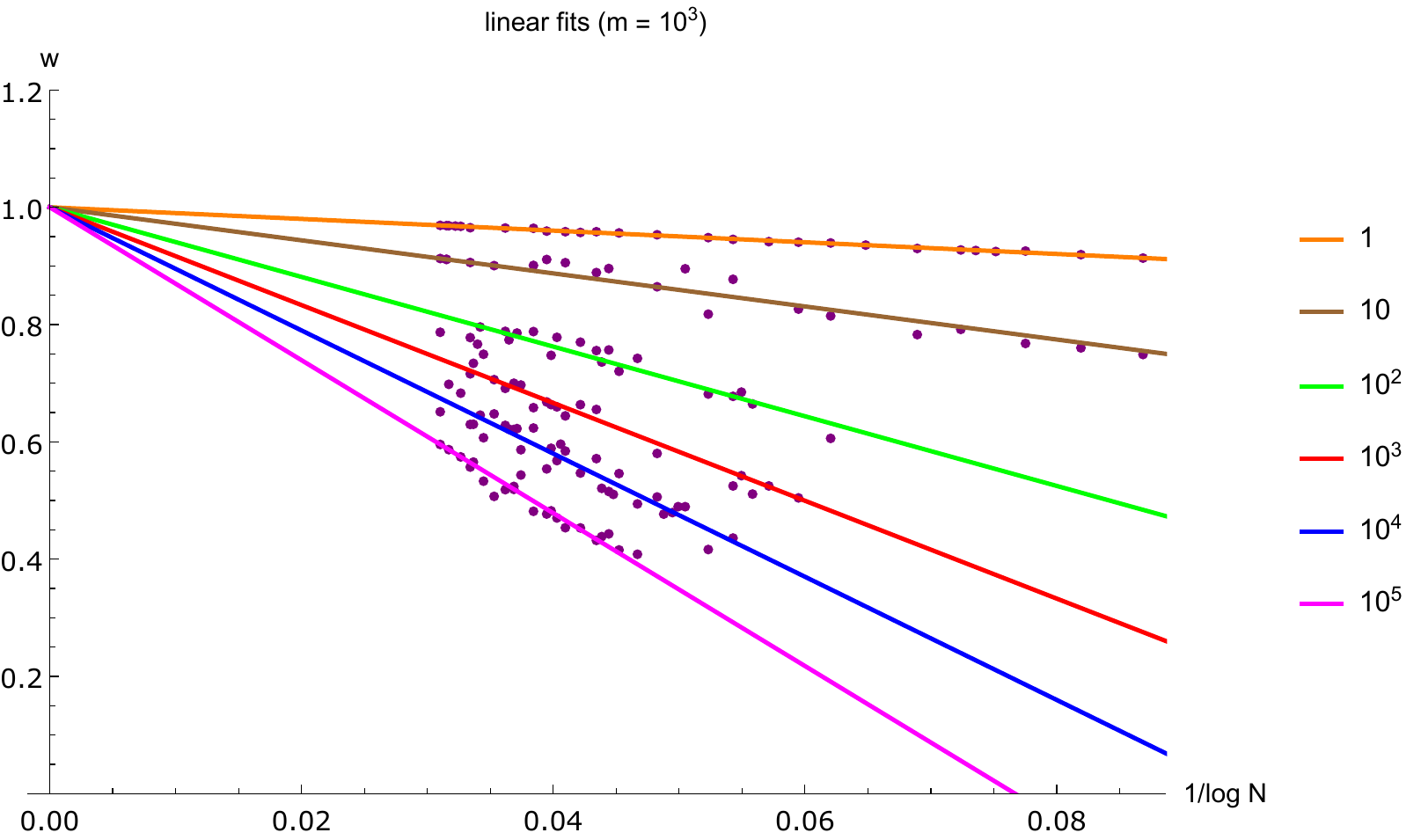}
                \caption{$m=10^{3}\; fits$}
                \label{fig:nperiodic1}
            \end{subfigure}\hfill
             \begin{subfigure}[t]{0.45\textwidth}
         \includegraphics[width=\textwidth]{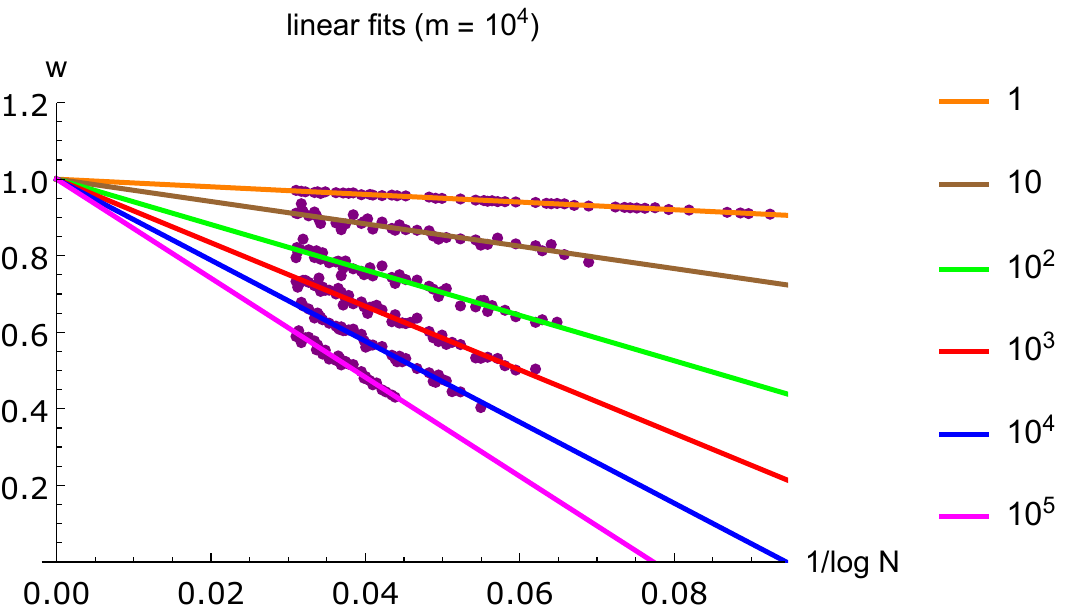}
                \caption{$m=10^{4}\; fits$}
                \label{fig:nperiodic2}
            \end{subfigure}
             \begin{subfigure}[t]{0.45\textwidth}
                \includegraphics[width=\textwidth]{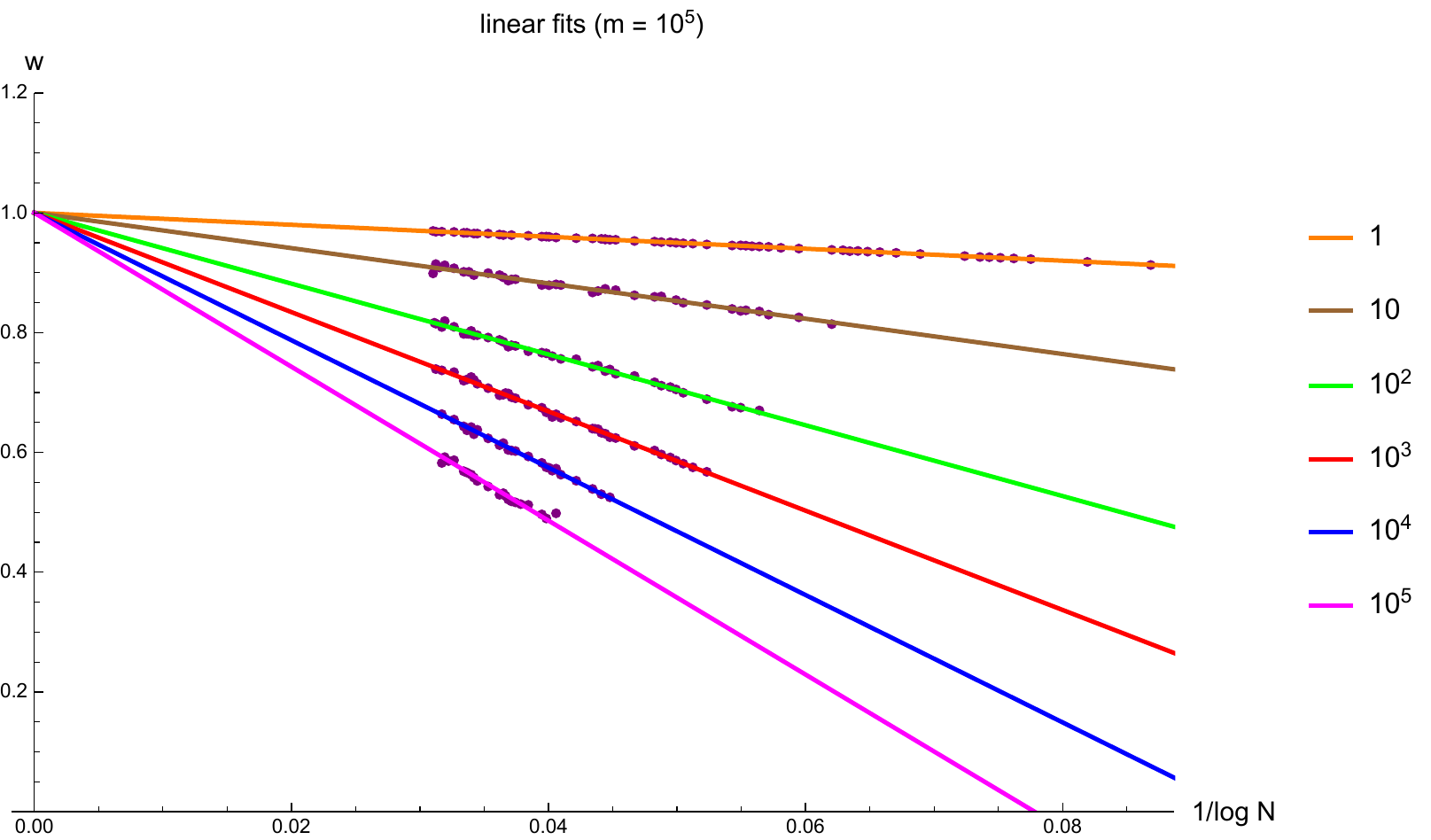}
                \caption{$m=10^{5}\; fits$}
                \label{fig:nperiodic3}
            \end{subfigure}
            \captionsetup{width=0.95\linewidth}
            \captionsetup{justification=centering}
              
\caption{Fits of $w$ for different h and a) $m=10^{3}$ b) $m=10^{4}$ c) $m=10^{3}$}
\label{fig:3}
\end{figure}

\section{Results}

\hspace{20.0pt} A comparison between the theoretical asymptotic and empirical values for the mean and variance of a sample interval is shown in Table 1. See that the systematic error has been neglected in the last columns.

\indent It is clear that the MS conjecture predicts the correct values of the mean, therefore no further analysis is needed. Nonetheless, it is worth mentioning that the statistical error is overestimated, this could be a consequence of using big amount of points to compute each parameter.
\begin{table}[!t]
\resizebox{\textwidth}{!}{%
\begin{tabular}{|l|l|l|l|l|l|l|l|}
\hline
b(h,m)     & m=$10^{3}$         & m=$10^{4}$         & m=$10^{5}$          & $\alpha (h,m)$ & m=$10^{3}$        & m=$10^{4}$        & m=$10^{5}$        \\ \hline
h=1        & 0.996 $\pm$ 0.004  & 0.998 $\pm$ 0.003  & 0.999 $\pm$ 0.001   & h=1            & -                 & -                 & -                 \\ \hline
h=$10$     & 2.819 $\pm$ 0.069  & 2.919 $\pm$ 0.037  & 2.946 $\pm$ 0.012   & h=$10$         & 0.846 $\pm$ 0.030 & 0.833 $\pm$ 0.016 & 0.845 $\pm$ 0.005 \\ \hline
h=$10^{2}$ & 5.941 $\pm$ 0.079  & 5.929 $\pm$ 0.041  & 5.911 $\pm$ 0.011   & h=$10^{2}$     & 1.071 $\pm$ 0.017 & 1.070 $\pm$ 0.009 & 1.071 $\pm$ 0.002 \\ \hline
h=$10^{3}$ & 8.343 $\pm$ 0.093  & 8.305 $\pm$ 0.035  & 8.288 $\pm$ 0.011   & h=$10^{3}$     & 1.062 $\pm$ 0.013 & 1.057 $\pm$ 0.005 & 1.055 $\pm$ 0.002 \\ \hline
h=$10^{4}$ & 10.503 $\pm$ 0.081 & 10.581 $\pm$ 0.029 & 10.635 $\pm$ 0.016  & h=$10^{4}$     & 1.031 $\pm$ 0.009 & 1.037 $\pm$ 0.003 & 1.046 $\pm$ 0.002 \\ \hline
h=$10^{5}$ & 13.041 $\pm$ 0.076 & 12.937 $\pm$ 0.036 & 12.8487 $\pm$ 0.037 & h=$10^{5}$     & 1.045 $\pm$ 0.007 & 1.036 $\pm$ 0.003 & 1.029 $\pm$ 0.003 \\ \hline
\end{tabular}}
\caption{Fits for b(h,m) and $\alpha(h,m)$}
\label{tab:2}
\end{table}

 \begin{figure}[!t]\centering
        \setlength\belowcaptionskip{-1.5ex}
          \begin{subfigure}[t]{0.45\textwidth}
           \includegraphics[width=\textwidth]{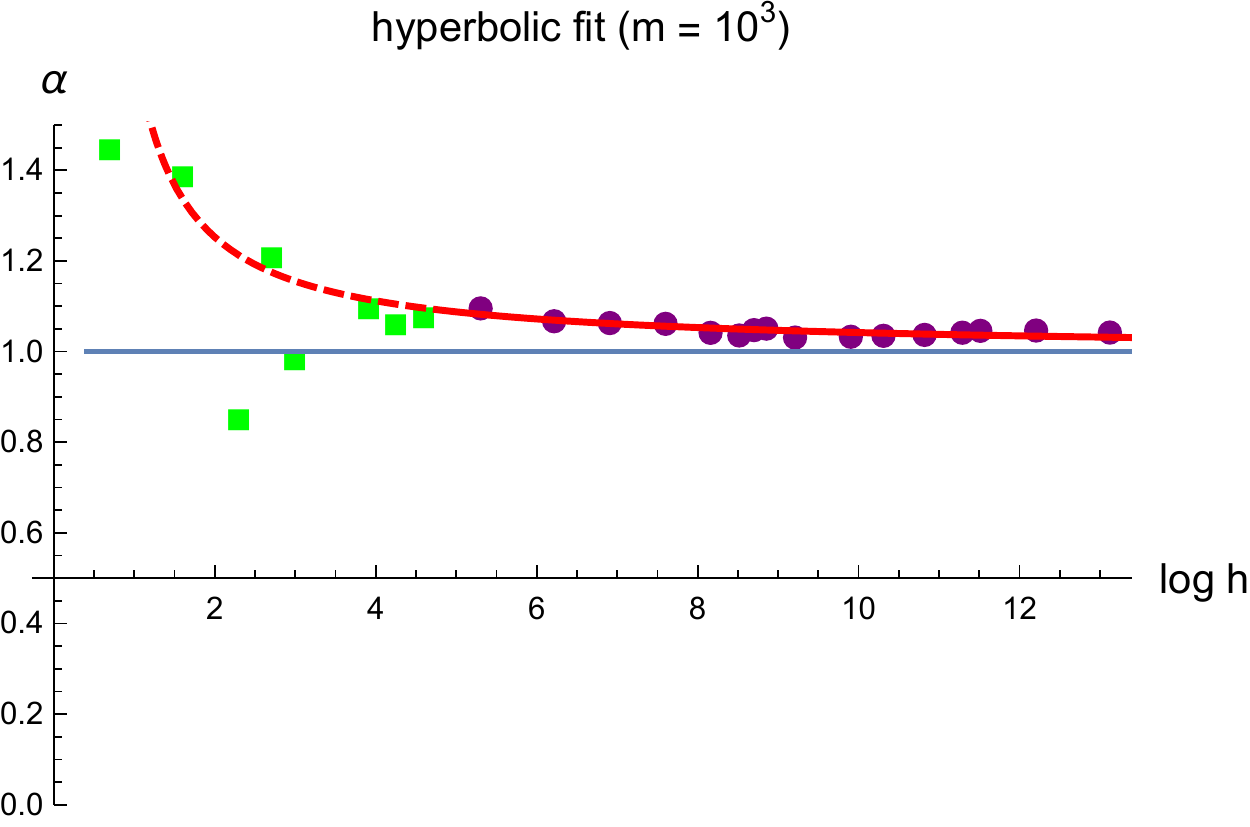}
                \caption{$\alpha(\log h)\; with \;m=10^{3}$}
                \label{fig:nperiodic1}
            \end{subfigure}\hfill
             \begin{subfigure}[t]{0.45\textwidth}
         \includegraphics[width=\textwidth]{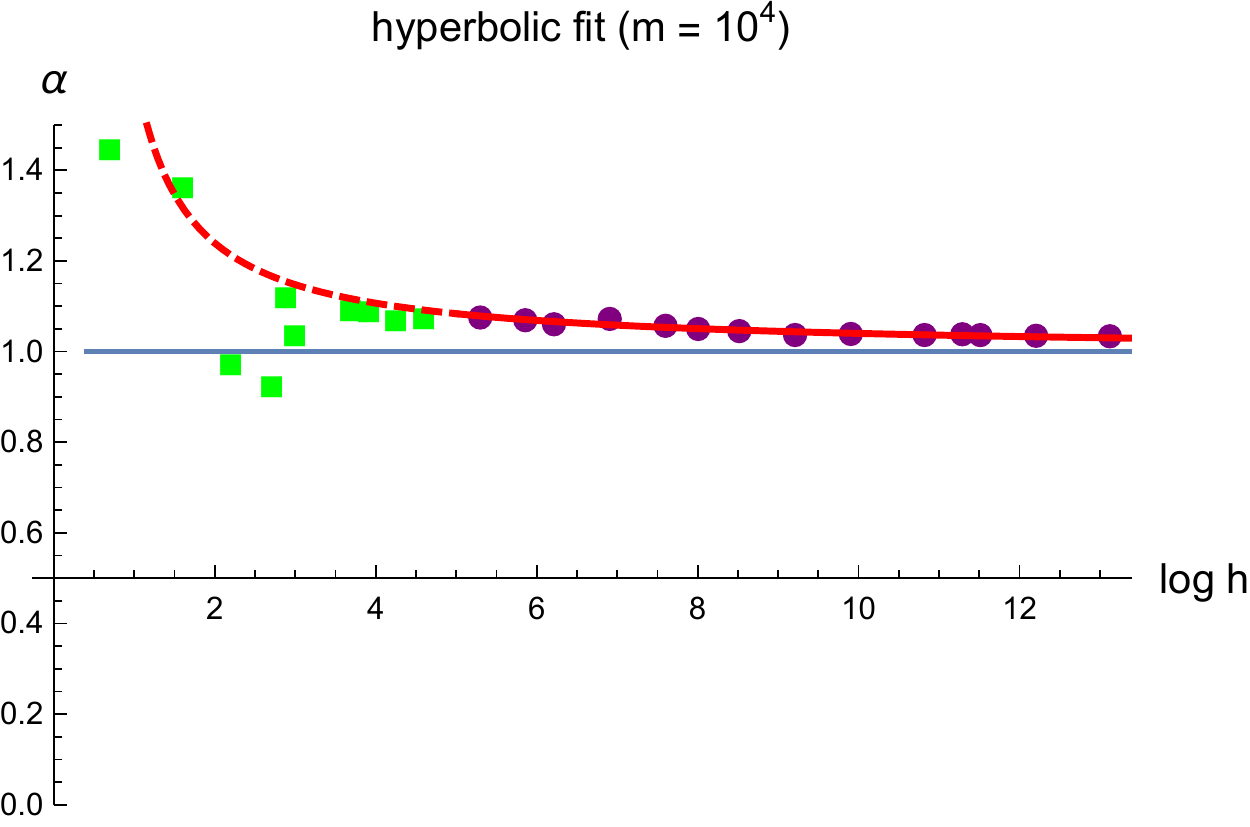}
                \caption{$\alpha(\log h)\; with\; m=10^{4}$}
                \label{fig:nperiodic2}
            \end{subfigure}
             \begin{subfigure}[t]{0.45\textwidth}
                \includegraphics[width=\textwidth]{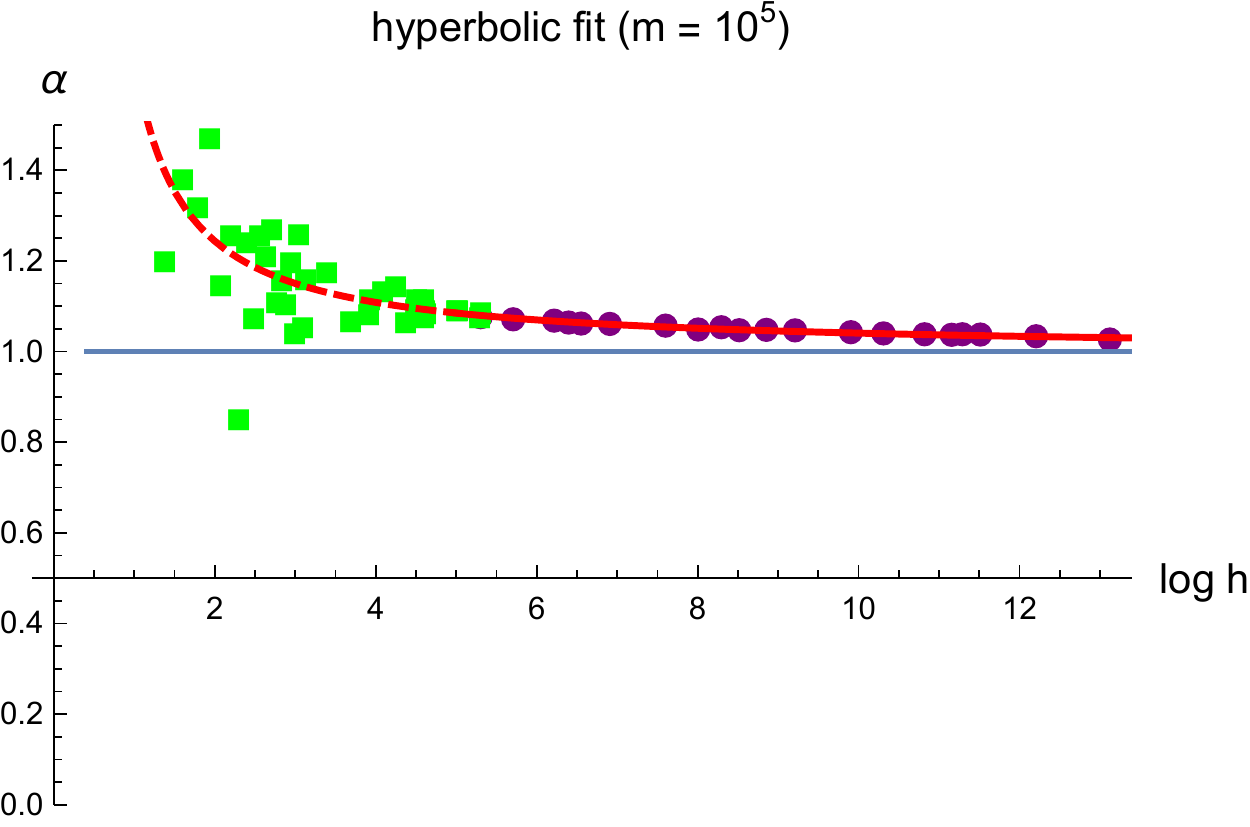}
                \caption{$\alpha(\log h)\; with\; m=10^{5}$}
                \label{fig:nperiodic3}
            \end{subfigure}
            \captionsetup{width=0.95\linewidth}
            \captionsetup{justification=centering}
              
\caption{Fits of $\alpha(\log h)$ for a) $m=10^{3}$ b) $m=10^{4}$ c) $m=10^{3}$}
\label{fig:4}
\end{figure}
\indent The results for the normalized variance can be seen in Figure \ref{fig:3}. The linear behaviour of $w$ is manifest, validating equation (\ref{eq:norvarg}). From the fitted values of $b(h,m)$ we can infer C and $\alpha(h)$ from the equation (\ref{eq:13}). All the values of $b(1,m)$ are coherent with 1, thence $C=1$. From Table \ref{tab:2} and Figure \ref{fig:4} we can see, as we expected, that $\underset{h \rightarrow  \infty}{\lim} \alpha (h)=1 $.  A simple way to achieve this result is by parametrizing $\alpha (h)=\frac{\alpha_{1}+f(h)}{\alpha_{2}+f(h)}$. See that there are oscillations for $h < 200$, nevertheless, these values are not included in the fit since they do not belong to the mesoscopic scale, since, on such a scale, $\frac{h}{\log N} \approx \infty$. The error bars do not appear in Figure \ref{fig:4} because they are negligible for $h>200$, as can be seen in Table \ref{tab:2}. We have to bear in mind that the dependence in m is just a matter of precision, but it can not explicitly appear in the formulas.

\indent From Table \ref{tab:2} we expect $\alpha (h)$ to be a monotonically decreasing function, with $f(h)$ being a monotonically increasing function. To verify both conditions it is necessary that $\alpha_{1}>\alpha_{2}$. We can rewrite $\alpha(h)$  as $\alpha (h)=\frac{\alpha_{1}+(f(h)^{-1})^{-1}}{\alpha_{2}+(f(h)^{-1})^{-1}}$, doing a Taylor expansion in $f(h)^{-1}$ around $h \rightarrow \infty$ we get:
\begin{equation}
   \alpha(h)=1+\frac{\alpha_{1}-\alpha_{2}}{f(h)}+ \mathcal{O}(\frac{1}{(f(h))^{2}})
   \label{eq:tay}
\end{equation}

\indent Comparing Equation (\ref{eq:tay}) with Equation (\ref{eq:13}) we obtain a variance $\sigma_{p} \approx\\\approx\frac{h}{\log^{2} N}(\log( \nicefrac{N}{\left( h^{\nicefrac{(\alpha_{1}-\alpha_{2})}{f(h)}+1}\right )})-1)$. Taking into account that $\sigma_{p} \sim \frac{h}{\log^{2} N} \log \frac{N}{h}$ it is reasonable to expect $f(h)=\log(h)$, as $\log\left (h^{\nicefrac{c}{\log(h)}} \right )=c$. In addition, logarithms are closely related with prime numbers \cite{10}. 

Apart from $\alpha (h)=\frac{\alpha_{1}+\log{h}}{\alpha_{2}+\log{h}}$, in appendix \ref{appendix:A} we have studied the parametrizations $\alpha_{\rm I} (h)=\frac{\alpha_{1}+\log{h}}{\log{h}}$ and $\alpha_{\rm II} (h)=\frac{\log{h}}{\alpha_{2}+\log{h}}$. Assessing the goodness of fit of each parametrization we have concluded that the best-fitting one is $\alpha_{\rm I} (h)=\frac{\alpha_{1}+\log{h}}{\log{h}}$. See Apendix \ref{appendix:A} for a detailed statistical analysis and Table (\ref{tab:3}) for the results of the fit. Identifying $B=1+\alpha_{1}-\alpha_{2}=1.414 \pm 0.004$ and substituting the expression (\ref{eq:tay}) in (\ref{eq:13}) we obtain:

%% \indent The most simple guess of $f(h)$ would be $f(h)=h$, nevertheless choosing that parametrization would imply working with a set of data formed by dozens of points in a range of the h-axis of 50.000 numbers and a range of the $\alpha$-axis of 2 numbers. That would lead to a rather unfriendly fit. \textcolor{blue}{As can be seen in the Appendix C} In fact, using the parametrization $f(h)=h$ results in a value for a and b of the order of $10^{6}$, with confidence intervals of the order of $10^{5}$. We are not stating that it is an erroneous parametrization but unsuitable for the amount of data that can be generated by a personal computer. \indent A convenient guess for $f(h)$ may be $f(h)=\log h$, as it solves the problem of the wide range of h ($\log 50.000 \approx 10.8$) 

\begin{equation}
   w=\frac{\log \nicefrac{N}{h}}{\log N}+\frac{B}{\log N}
   \label{eq:finall}
\end{equation}

\indent The variance can be obtained by multiplying equation (\ref{eq:finall}) by the mean ($\nicefrac{h}{\log N}$), which, in the asymptotic limit, tends to:

\begin{equation}
   \underset{h \rightarrow  \infty}{\lim} \sigma^{2}=\frac{h}{\log^{2}N}\log \frac{N}{h}
   \label{eq:finalLim}
\end{equation}

\indent Which is the result predicted by the MS conjecture.

\begin{table}[t]
\centering
\begin{tabular}{|c|c|c|c|}
\hline
m & $10^{3}$           & $10^{4}$          & $10^{5}$          \\ \hline
$B$ & $1.420 \pm 0.023$  & $1.408 \pm 0.010$ & $1.414 \pm 0.004$ \\ \hline
\end{tabular}
\captionsetup{justification=centering,margin=2cm}
\caption{Fits for $B=1+\alpha_{1}$. Using\\ $\alpha_{\rm I}(h)=\frac{\alpha_{1}+\log{h}}{\log{h}}$ }
\label{tab:3}
\end{table}

\section{Conclusion}

\hspace{20.0pt} This project aimed to verify experimentally the MS conjecture and to provide an empirical formula to extend its range of application on the mesoscopic scale. We have obtained a formula that predicts the variance of the distribution of prime numbers on such scale for large but finite N.

\begin{equation}
    \sigma^{2}=\frac{h}{\log^{2}N}(\log \frac{N}{h}+B)
   \label{eq:finalLim}
\end{equation}

\indent This formula is completely consistent with the asymptotic limit of the MS conjecture. We have computed $B=1.414 \pm 0.004$ using the parameters of the best-fitting model. It is noteworthy to mention that, in all parametrizations, $B$ is coherent with $-B'=\gamma+\log{2\pi}-1 \approx 1.415$, where $\gamma$ is the Euler–Mascheroni constant, which arises in the calculation of the moments of the distribution of prime numbers in some scales \cite{7}. See that we have obtained the same result as Sanchis \cite{9} using a different parametrization of $\alpha(h)$.

\indent All the calculations have been done using \emph{Mathematica} and have been limited by the computation capacity of an average computer. For that reason, we have worked with N less than $10^{14}$. In order to verify the range of validity of the equation (\ref{eq:finalLim}), it would be necessary to carry out calculations with bigger N, that would require a larger computational capacity. Nevertheless, in the range in which we have been able to work, we have provided an experimental result that supports the validity of the MS conjecture.

\clearpage

\begin{appendix}

\section{Statistical analysis of the parametrizations $\alpha(h)$}
\label{appendix:A}

\hspace{20.0pt} In the parametrization of $\alpha(h)$ we are only interested in the behaviour at large h, hence the only significant parameter, according to equation (\ref{eq:tay}), will be $B=1+\alpha_{1}-\alpha_{2}$. To assess the goodness of fit of every parameterization we are going to employ the Pearson's $\chi^{2}$ test to compute the p-value \cite{11}. The p-value would give us an estimation of the degree of confidence in the null hypothesis, in this case, that hypothesis corresponds to consider that $\alpha(h)$ is well fitted by the parametrization. We have used 16 points for $m=10^{3}$, 14 points for $m=10^{4}$ and 20 points for $m=10^{5}$. Let $\nu$ be the degrees of freedom of the model, that may be computed as the number of points minus the number of parameters. In the tables of this appendix we have shown values of the different parameters used in the different parametrizations and the reduced chi-square ($\chi_{\nu}^{2}=\nicefrac{\chi^{2}}{\nu}$) and p-value for each one: 

\begin{table}[h]
\centering
\begin{tabular}{|c|c|c|c|}
\hline
m & $10^{3}$           & $10^{4}$          & $10^{5}$          \\ \hline
$\alpha_{1}$ & $-0,475 \pm 1.364$  & $0.571 \pm 0.821$ & $0.666 \pm 0.305$ \\ \hline
$\alpha_{2}$ & $-0.848 \pm 1.291$ & $0.154 \pm 0.776$ & $0.238 \pm 0.289$ \\ \hline
$B$ & $1.373 \pm 1.878$  & $1.597 \pm 1.129$ & $1.428 \pm 0.420$ \\ \hline
\end{tabular}
\captionsetup{justification=centering,margin=2cm}
\caption{Fits for $B=1+\alpha_{1}-\alpha_{2}$. Using\\ $\alpha_{\rm III}(h)=\frac{\alpha_{1}+\log{h}}{\alpha_{2}+\log{h}}$ }
\label{tab:as}
\end{table}
\begin{table}[h]
\centering
\begin{tabular}{|c|c|c|c|}
\hline
m & $10^{3}$           & $10^{4}$          & $10^{5}$          \\ \hline
$B$ & $1.420 \pm 0.023$  & $1.408 \pm 0.010$ & $1.414 \pm 0.004$ \\ \hline
\end{tabular}
\captionsetup{justification=centering,margin=2cm}
\caption{Fits for $B=1+\alpha_{1}$. Using\\ $\alpha_{\rm I}(h)=\frac{\alpha_{1}+\log{h}}{\log{h}}$ }
\label{tab:asd}
\end{table}
\begin{table}[!h]
\centering
\begin{tabular}{|c|c|c|c|}
\hline
m & $10^{3}$           & $10^{4}$          & $10^{5}$          \\ \hline
$B$ & $1.399 \pm 0.021$  & $1.386 \pm 0.009$ & $1.392 \pm 0.004$ \\ \hline
\end{tabular}
\captionsetup{justification=centering,margin=2cm}
\caption{Fits for $B=1-\alpha_{2}$. Using\\ $\alpha_{\rm II} (h)=\frac{\log{h}}{\alpha_{2}+\log{h}}$ }
\label{tab:sa}
\end{table}
\begin{table}[!h]
\centering
\begin{tabular}{|c|c|c|c|c|c|c|c|}
\hline
$\chi_{\nu}^{2}$     & m=$10^{3}$         & m=$10^{4}$         & m=$10^{5}$          & p-value & m=$10^{3}$        & m=$10^{4}$        & m=$10^{5}$        \\ \hline
$\alpha_{\rm I}$        & 1.683 & 1.400  & 1.063  & $\alpha_{\rm I}$            & 4.6 \%                 & 15 \%                 & 38\%                 \\ \hline
$\alpha_{\rm II}$      & 1.718 & 1.447 & 1.351  & $\alpha_{\rm II}$         & 4 \% & 13\% & 14\% \\ \hline
$\alpha_{\rm III}$  & 1.904  & 1.521 & 1.078   & $\alpha_{\rm III}$      & 2 \% & 11 \% & 37 \% \\ \hline
\end{tabular}
\captionsetup{justification=centering,margin=2cm}
\caption{$\chi_{\nu}^{2}$ and p-value for all the parametrizations}
\label{tab:def}
\end{table}

\indent See that, for $\alpha_{\rm III}(h)$, there is no contradiction in the fact of having different values of  $\alpha_{1}$ and $\alpha_{2}$ for different m. As stated before, we are only interested in the difference $\alpha_{1}-\alpha_{2}$, that is coherent between the different sets characterized by $m$ in each of the three parametrizations. From Table (\ref{tab:def}) we can see that the best-fitting parametrization is $\alpha_{\rm I}(h)=\frac{\alpha_{1}+\log{h}}{\log{h}}$. We can infer the value $B=1.414 \pm 0.004$.

\end{appendix}
%Comenta lo de los limites de mathematica.
%Comenta que consideramos un fitting suitable if it reproduces the MS conjecture in assympt limit
%Teoría: hemos modelado la varianza normalizada porque es más sencilla que la farianza (ya que no salen cosas al cuadrado)
 
\end{document}